\documentclass[12pt]{article}
\usepackage[mathscr]{eucal}
\usepackage{amsmath,amssymb,amscd}
\usepackage{color}
\usepackage{epsfig}

\def \beq{ \begin{equation}}
\def \eeq{\end{equation}}

\begin{document}

\begin{center}
{\bf A dissection proof of the law of cosines,\\
replacing
Cuoco-McConnell's rectangles with congruent triangles.}
\end{center}

\begin{center}
{\bf Martin Celli}
\end{center}

\begin{center}
March 31st 2018
\end{center}

\begin{center}
Departamento de Matem\'aticas\\
Universidad Aut\'onoma Metropolitana-Iztapalapa\\
Av. San Rafael Atlixco, 186. Col. Vicentina. Del. Iztapalapa. CP
09340. Mexico City.\\
E-mail: cell@xanum.uam.mx
\end{center}

\vskip1cm

Taking up the challenge McConnell laid down at the end of his proof
of the law of cosines ([6]), we give a completely visual dissection
proof of this theorem, which applies to any triangle. In order to
avoid the tri\-go\-no\-me\-tric expressions of Cuoco-McConnell's proof, we
replaced the equal-area rectangles with congruent triangles. As a
matter of fact, tri\-go\-no\-me\-tric expressions are implicitly based on
the similarity of two right triangles with a common non-right angle.
So they are conceptually less simple than our congruent triangles
which are, moreover, easy to visualize. This makes our proof the
only dissection proof and the simplest proof of its family, and thus
one of the best options for a course of geometry.\\

Cuoco-McConnell's proof was published in [6], [3], [5], [8]. In [8],
the proof is written in an algebraic way, whereas in [6], [3], [5],
it is written in an equivalent geometric way involving three pairs
of rectangles with the same area. A similar proof, involving two
pairs of rectangles with the same area, can be found in [7].
Unfortunately, in these references, the areas of all the rectangles
had to be expressed as tri\-go\-no\-me\-tric functions, in order to show
that they were pairwise equal. In the particular case of a right
triangle (where a pair of rectangles have zero area), and as in our
proof, Euclid had already managed to avoid this by considering, for
every pair of rectangles, a pair of congruent triangles with half
the area of one of the rectangles ([4]). A generalization of
Euclid's proof to the law of cosines can be found in [2] (where
every triangle is replaced with a parallelogram made of this
triangle and its symmetric). Each triangle of Euclid, as well as
each triangle of our proof, has two sides in common with the
triangle \(ABC\) of the theorem. However, denoting by \(x\) the
angle formed by these sides, the corresponding angle in Euclid's
triangle is \(\pi /2+x\), whereas the corresponding oriented angle
in our triangle is \(\pi /2-x\). As a consequence of this, and
unlike Euclid's triangles and the parallelograms of [2], two
triangles of a same pair of our proof do not have intersection. This
helps us in finding a visual dissection argument (cutting each
rectangle into a set of congruent pieces) in order to prove that the
area of the triangle is half the area of the corresponding
rectangle. In both our proof and [2], tri\-go\-no\-me\-try is only used to
compute the inevitable area of the triangles/parallelograms
corresponding to the term \(2CA.CB\cos (BCA)\) of the identity.\\

When one of the points \(A'\), \(B'\), \(C'\) of our proof
(symmetric of \(A\), \(B\), \(C\) with respect to \(BC\), \(CA\), \(AB\))
is located outside its square, some areas need to be interpreted in an oriented sense,
and counted as negative. This already happened in [6], [3], [5],
when the triangle was obtuse. Here, it can also happen when an altitude
of the triangle \(ABC\) is longer than its perpendicular side.
We give the proof in the acute case, when the altitude passing through \(B\)
is longer than \(AC\) (thus, the point \(B'\) is outside its square),
the other cases can be studied in a similar way.\\

It is worth mentioning that, in the particular case where the three
points \(A'\), \(B'\), \(C'\) are located inside their squares, our
proof becomes equivalent to Anderson's ([1]), where both small
triangles adjacent to a same blue/green/red triangle of the figure
are cut and pasted, so that the three triangles (or equivalently,
the blue/green/red triangle and its symmetric) make a parallelogram.
Thus, applying the same transformation (replacing a triangle with a
parallelogram made of it and its symmetric), we can obtain the
parallelograms of [2] from Euclid's triangles, and Anderson's
parallelograms from the triangles of the proof of this note. As
every triangle of this proof has two sides in common with a triangle
of Euclid, with a corresponding angle \(\pi-\alpha =\pi /2-x\)
instead of \(\alpha =\pi /2+x\), every parallelogram of Anderson
and its corres\-ponding parallelogram of [2] are congruent.\\

\begin{figure}[h]
  \hspace{1.5cm}
   \includegraphics[width=25cm]{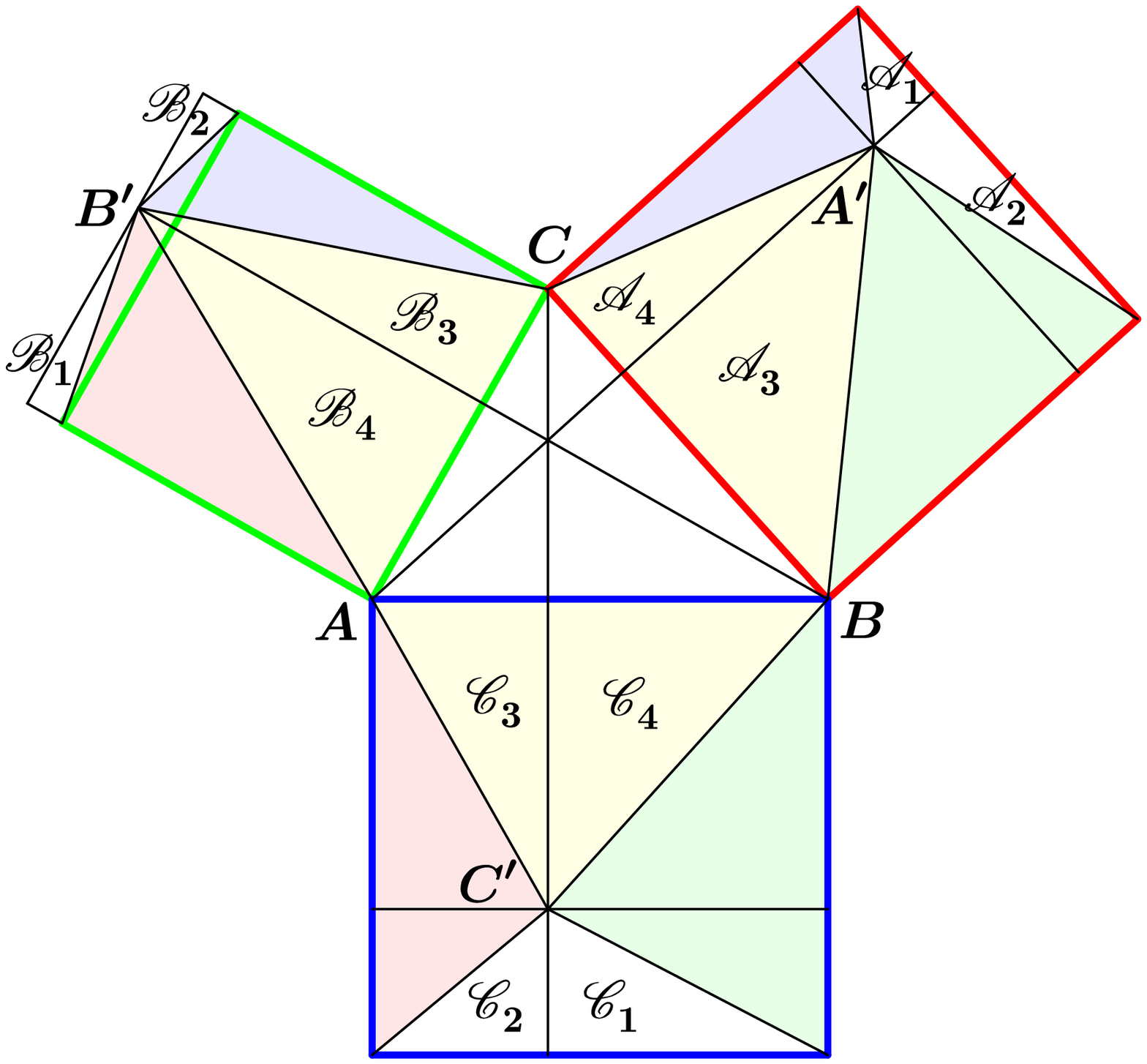}
    \end{figure}

Here is our proof:
\[BC^2+CA^2-AB^2=[\mbox{red square}]+[\mbox{green square}]-[\mbox{blue square}]\]
\[=(2\mathscr{A}_3+2\mathscr{A}_2+2\mathscr{A}_4+2\mathscr{A}_1)
+((2\mathscr{B}_3-2\mathscr{B}_2)+(2\mathscr{B}_4-2\mathscr{B}_1))
-(2\mathscr{C}_3+2\mathscr{C}_2+2\mathscr{C}_4+2\mathscr{C}_1)\]
\[=2((\mathscr{A}_3+\mathscr{A}_2)+(\mathscr{A}_4+\mathscr{A}_1)
+(\mathscr{B}_3-\mathscr{B}_2)+(\mathscr{B}_4-\mathscr{B}_1)
-(\mathscr{C}_3+\mathscr{C}_2)-(\mathscr{C}_4+\mathscr{C}_1))\]
\[=2([\mbox{green triangle}]+[\mbox{blue triangle}]+[\mbox{blue triangle}]
+[\mbox{red triangle}]-[\mbox{red triangle}]-[\mbox{green triangle}])\]
\[=4[\mbox{blue triangle}]=2CA.CB\sin (\pi /2-BCA)=2CA.CB\cos(BCA)\cdot \]

\(\)\\

{\bf References.}\\
\(\mbox{[1]}\) S. Anderson et al. Comments about McConnell's
"Illustrated Law of Cosines".
Website "Interactive Mathematics Miscellany and Puzzles" (A. Bogomolny):\\
https://www.cut-the-knot.org/htdocs/dcforum/DCForumID11/90.shtml\\
\(\mbox{[2]}\) I. Boyadzhiev. The Law of Cosines.
Website "Discover Math with GeoGebra":\\
https://www.geogebra.org/m/ScsxfNx3\\
\(\mbox{[3]}\) A. Cuoco, A. Baccaglini-Frank, J. Benson, N. Antonellis D'Amato, D. Erman,
B. Harvey, K. Waterman. CME Project: Geometry. Pearson Education (2009).\\
\(\mbox{[4]}\) Euclid. Elements. Book I, proposition 47. Can be found in:\\
https://mathcs.clarku.edu/\~{ }djoyce/java/elements/bookI/propI47.html\\
\(\mbox{[5]}\) R. E. Howe. The Cuoco Configuration.
Amer. Math. Monthly 120 (2013), no. 10, 916-923.\\
\(\mbox{[6]}\) D. McConnell. The Illustrated Law of Cosines.
Website "Interactive Mathematics Miscellany and Puzzles" (A. Bogomolny):\\
https://www.cut-the-knot.org/pythagoras/DonMcConnell.shtml\\
\(\mbox{[7]}\) J. Molokach. Law of Cosines-A Proof Without Words.
Amer. Math. Monthly 121 (2014), no. 8, 722.\\
\(\mbox{[8]}\) Entry "Law of cosines" (3.2 Proofs, using trigonometry).
Website "Wikipedia":\\
https://en.wikipedia.org/wiki/Law\_of\_cosines\#Using\_trigonometry\\

\(\)\\

Martin Celli.\\
Depto. de Matem\'aticas, Universidad Aut\'onoma
Metropolitana-Iztapalapa. Mexico City.\\
E-mail: cell@xanum.uam.mx

\end{document}